\documentclass[final,5p,times,twocolumn]{elsarticle}

\usepackage{graphicx} 
\usepackage{mathtools}
\usepackage{amsmath, nccmath}
\usepackage{mathrsfs}
\usepackage{xcolor}
\usepackage{comment}
\usepackage{amssymb}
\usepackage{float}
\usepackage{subfigure}
\usepackage{epstopdf}
\usepackage{lineno,hyperref}
\usepackage{color}
\usepackage{tikz,amsmath,siunitx}
\usetikzlibrary{arrows,backgrounds,patterns,matrix,shapes,fit,calc,shadows,plotmarks}

\newtheorem{thm}{Theorem}

\definecolor{ForestGreen}{rgb}{0.3, 0.7, 0.3}

\modulolinenumbers[5]

\journal{arxiv.org}

\begin{document}

\begin{frontmatter}

\title{\textbf{Extremum Seeking for Stefan PDE with Moving Boundary}}


\author[uerj]{Maurício Linhares Galvão}
\ead{mauriciolinhares92@gmail.com}

\author[uerj]{Tiago Roux Oliveira}
\address[uerj]{Dept. of Electronics and Telecommunication Engineering, State University of Rio de Janeiro (UERJ), Rio de Janeiro, RJ 20550-900, Brazil.}
\ead{tiagoroux@uerj.br}

\author[ucsd]{Miroslav Krstic}
\address[ucsd]{Dept. of Mechanical and Aerospace Engineering, University of California - San Diego (UCSD), La Jolla, CA 92093-0411, USA.}
\ead{krstic@ucsd.edu}

%
%


\begin{abstract}
This paper presents the design and analysis of the extremum seeking for static maps with input passed through a partial differential equation (PDE) of the diffusion type defined on a time-varying spatial domain whose boundary position is governed by an ordinary differential equation (ODE). This is the first effort to pursue an extension of extremum seeking from the heat PDE to the Stefan PDE. We compensate the average-based actuation dynamics by a controller via backstepping transformation for the moving boundary, which is utilized to transform the original coupled PDE-ODE into a target system whose exponential stability of the average equilibrium of the average system is proved. The discussion for the delay-compensated extremum seeking control of the Stefan problem is also presented and illustrated with numerical simulations. 
\end{abstract}



\end{frontmatter}


\section{Introduction}
Extremum seeking (ES) is a non-model based approach in the field of adaptive control which searches in real time extremum points (maximum or minimum) of a performance index of a system. This method has received great attention in the control community by the means of facing control problems when the plant has imperfections in its model or uncertainties dynamics \citep{Krsti2000StabilityOE}.

ES was first introduced in \citep{leblanc1922} for maximizing power transfer to a tram car. Along the history, the number of publications concerning ES remained low until it's first general stability proof for stable dynamic systems with unknown output maps was carried out \citep{Krsti2000StabilityOE}. Since then, important studies in theory and applications were developed, such as  \citep{transient_SCL,Ariyur2003RealTimeOB,source_seeking,Ghaffari2011MultivariableNE,Manzie2009ExtremumSW,Scheinker_TAC_2013,SK:2014,ghaffari_TCST_2014,paz_2020,JPC2,JPC1}.

Reference \citep{Oliveira2017ExtremumSF} was the first one to handle partial differential equations (PDEs) with ES scheme, addressing the design and analysis of multi-variable static maps subject to arbitrarily long time delays. The delays pointed out by the authors can be modelled as first-order hyperbolic transport PDEs \cite{optimality}. This idea has enabled the development of extension to other classes of PDEs with fixed domains \cite{book2022,cascades}.

On the other hand, a large number of applications in various areas appear as moving boundary or phase change problems, such as \citep{wettlaufer1991heat} and \citep{petrus2012enthalpy}. Usually, these kind of problems arise in heat conduction situations and need to be solved in a time-dependent space domain with a moving boundary condition. For this reason diffusion PDEs with moving boundaries, known as ``Stefan problem" have been studied actively for the last few decades \cite{JPC3}, \cite{JPC4}.

The dynamics of the the position of the moving boundary in the Stefan problem  is governed by an ODE which depends on the PDE state, generating a nonlinear coupling of  the PDE and ODE dynamics, increasing the complexity of the problem when compared to conventional analyses for PDEs of fixed domains (not depending on time or states) and ODEs.

In this paper, we develop an ES controller for the Stefan problem. The integrator that is usually employed in the ES scheme can be leveraged as part of the Stefan model, just like the one proposed in \citep{zhang2007extremum}. The objective of the ES will be to find the maximizer interface $s^{*}$ of some unknown map $Q(s^{*})$ aiming to regulate the phase change interface position to a value that attains the extremum. For this purpose, we design a compensator of the heat PDE with moving boundary and the probing signal, which is the result of solving the problem of generating a sinusoid at the distal end of a boundary-actuated heat equation. As a further contribution, we also study the effect of delays in the PDE actuation dynamics by compensating it via predictor feedback \citep{bresch2010delay}. 
 
An important discussion is about the validation of the Stefan model. Although the usual sinusoidal movement provoked by the ES algorithm may violate the phase maintenance when the extremum is achieved or during the transient either, we can keep the phase maintenance at least for the average system, thus, preserving the convergence analysis. The main theoretical contribution with respect to the previous conference paper \cite{TDS_2022_galvao} is the inclusion of the proof for the main theorem not presented before due to space limitations.

%
%


\newpage
This paper is organized as follows. The problem statement is presented in Sections \ref{one-phase Stefan problem} and \ref{extrachapter.ESC}, containing the relation of the ES scheme with the one-phase Stefan problem. In Section \ref{extrachapter.controller}, the control related problem is designed. Section \ref{extrachapter.stability} describes the stability analysis of the closed-loop average system with exponential convergence in $\mathscr{H}_{1}$-norm of the distributed temperature. The convergence of the moving boundary to the desired equilibrium is discussed in Section~\ref{samsungs23+}. Section~\ref{delay_coreano} shows the design of a compensator of the transport-heat PDE cascade in case of delays in the actuation dynamics. Numerical simulations are provided in Section \ref{extrachapter.simulation}, followed by conclusions in Section \ref{concludingfinally}.

\smallskip

\textbf{Norms and Notations:} We denote the partial derivatives of a function $u(x, t)$ as $\partial_{x} u(x,t) = \partial u(x,t)/\partial x$, $\partial_{t} u(x,t) = \partial u(x,t)/\partial t$. We conveniently use the compact form $u_{x}(x,t)$ and $u_{t}(x,t)$ for the former and the latter, respectively. The two-norm of a finite-dimensional ODE state vector $X(t)$ is denoted by single bars, $|X(t)|$. In contrast, norms of functions (of $x$) are denoted by double bars. We denote the spatial $\mathcal{L}_{2}[0,D]$ norm of the PDE state $u(x, t)$ as $\|u(t)\|^{2}_{\mathcal{L}_{2}([0,D])} \coloneqq \int_{0}^{D} u^2(x,t) \,dx$, where we drop the index $\mathcal{L}_{2}([0,D])$ in the following, hence $\|\cdot\| = \|\cdot\|_{\mathcal{L}_{2}([0,D])}$, if not otherwise specified. Moreover, the $\mathscr{H}_{1}$-norm is given by $\|u(t)\|^{2}_{\mathscr{H}_{1}} = \|u(t)\|^{2}_{\mathcal{L}_{2}}+\|u_{x}(t)\|^{2}_{\mathcal{L}_{2}}$. As defined in \citep{KH:02}, a vector function $f(t,\epsilon) \in \mathbb{R}^{n}$ is said to be of order $\mathcal{O}(\epsilon)$ over an interval $[t_{1},t_{2}]$, if $\exists k, \bar{\epsilon}$ : $|f(t,\epsilon)| \leq k \epsilon, \forall t \in [t_{1},t_{2}]$. In most cases, we provide no precise estimates for the constants $k$ and $\bar{\epsilon}$, and we use $\mathcal{O}(\epsilon)$ to be interpreted as an order of magnitude relation for sufficiently small $\epsilon$.

\section{One-phase Stefan Problem\label{one-phase Stefan problem}}

\begin{figure}[t]
\begin{center}
\includegraphics[width=8.5cm]{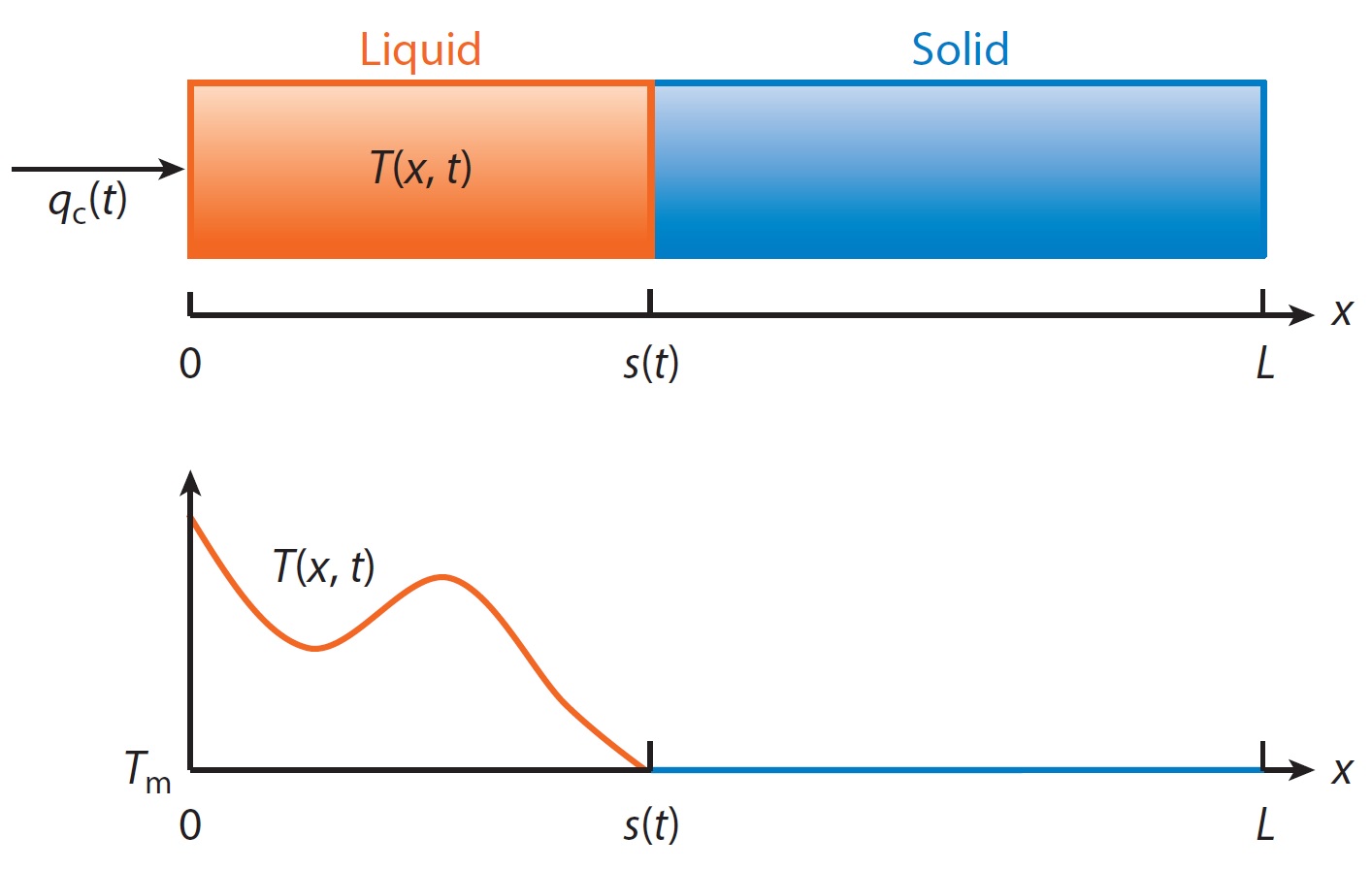}    
\caption{Schematic of one-phase Stefan problem \cite{Shumon_green_book:2020,AR_KK:2022}.
The temperature profile in the solid phase is assumed to be a uniform melting temperature.}
\label{extrachapter.fig:stefan}
\end{center}
\end{figure}

\begin{figure}[!htb]
    \begin{center}
    \includegraphics[width=8.5cm]{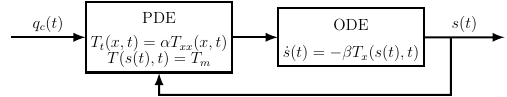} 
    \caption{The cascade of the PDE dynamics and the ODE system.}
    \label{extrachapter.fig:ode_plant}
    \end{center}
\end{figure}

The physical model which describes the 1-D Stefan problem in a pure one-component material of length $L$ is described in Figure~\ref{extrachapter.fig:stefan}. The domain $[0,L]$ is divided in two sub-domains $[0,s(t)]$ and $[s(t),L]$ which represents the liquid phase and the solid phase, respectively. The system is controlled by the heat flux $q_{c}(t)$ at $x=0$, because we are dealing with a Neumann boundary actuation as shown below:
\begin{align}
    T_{t}(x,t) &= \alpha T_{xx}(x,t),\quad x\in (0,s(t)),\quad \alpha = \dfrac{k}{\rho C_{p}}\label{extrachapter.T1}\\
    -kT_{x}(0,t) &= q_{c}(t)\label{extrachapter.T2}\\
    T(s(t),t) &= T_{m}\label{extrachapter.T3}\\
    \dot{s}(t) &= -\beta T_{x}(s(t),t),\quad \beta = \dfrac{k}{\rho \Delta H^{*}}, \label{extrachapter.T4}
\end{align}
where $T(x,t)$, $T_{m}$, $q_{c}(t)$, $k$, $\rho$, $C_{p}$ and $\Delta H^{*}$ are the distributed temperature of the liquid phase, melting temperature, manipulated heat flux, liquid heat conductivity, liquid density, liquid heat capacity and latent heat of fusion, respectively. Equations (\ref{extrachapter.T2}) and (\ref{extrachapter.T3}) are the boundary conditions of the system and (\ref{extrachapter.T4}) is the Stefan condition, which describes the dynamic of the moving boundary. Figure \ref{extrachapter.fig:ode_plant} shows the block diagram of the PDE-ODE cascade represented by equations (\ref{extrachapter.T1})-(\ref{extrachapter.T4}).

\section{Problem Statement} \label{extrachapter.ESC}

\begin{figure*}[!htb]
\begin{center}
\includegraphics[width=14.0cm]{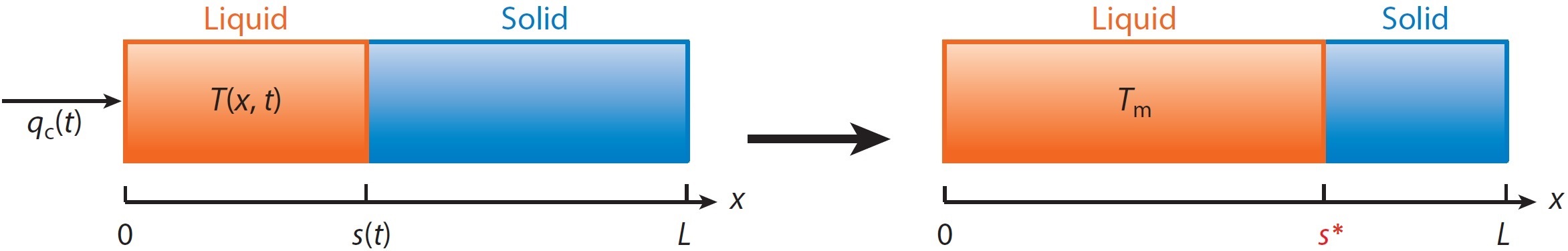}    
\caption{Control objective of the Stefan problem. We aim to design a heat flux input $q_c(t)$ such that the interface position $s(t)$ is driven to the setpoint position $s^*$.}
\label{extrachapter.fig:stefan_objectives}
\end{center}
\end{figure*}

For the sake of simplicity, we consider actuation dynamics which are described by a heat equation with $\alpha,\beta,k = 1$, $\theta(t) \in \mathbb{R}$ and the propagated actuator $\Theta(t) \in \mathbb{R}$ given by 
\begin{align}
    \dot{\Theta}(t) = \dot{s}(t) &= -\alpha_{x}(s(t),t)),\quad x \in (0,s(t))\label{extrachapter.Theta_p_actuator}\\
    \partial_{t}\alpha(x,t) &= \partial_{xx}\alpha(x,t)\label{extrachapter.stefan_actuator}\\
    \alpha(s(t),t) &= 0\label{extrachapter.boundary_actuator_0}\\
    -\partial_{x}\alpha(0,t) &= \theta(t)\label{extrachapter.theta_actuator},
\end{align}
where $\alpha:[0,s(t)]\times\mathbb{R}_{+}\rightarrow\mathbb{R}$ is $\alpha(x,t) = T(x,t) - T_{m}$ and $s(t) = \Theta(t)$ is the unknown interface represented as the moving boundary. The output is measured by the unknown static map with input (\ref{extrachapter.Theta_p_actuator}): 
\begin{equation}
y(t) = Q(\Theta(t)).
\label{extrachapter.eq:initial_output_static_map}
\end{equation}
The ES goal is to optimize an unknown static map Q($\cdot$) using a real-time optimization control with optimal unknown output $y^{*}$ and optimizer $\Theta^{*}$ as well as measurable output $y$ and input $\theta$. Consequently, the control objectives of the Stefan problem are achieved, \textit{i.e.}, $\lim\limits_{t\rightarrow \infty}s(t) = s^*$ and $\lim\limits_{t\rightarrow \infty}T(x,t) = T_m\,, \forall x \in [0, s^*]$, as illustrated in Figure~\ref{extrachapter.fig:stefan_objectives}. 

The unknown nonlinear map is locally quadratic, such that
\begin{equation}
Q(\Theta) = y^{*} + \dfrac{H}{2}(\Theta - \Theta^{*})^{2},
\label{extrachapter.eq:static_map}
\end{equation}
where $\Theta^{*},y^{*} \in \mathbb{R}$ and $H<0$ is the Hessian. Hence, the output of the static map is given by
\begin{equation}
y(t) = y^{*} + \dfrac{H}{2}(\Theta(t) - \Theta^{*})^{2}.
\label{extrachapter.eq:final_output_static_map}
\end{equation}


\subsection{Demodulation Signals}

The demodulation signal $N(t)$ which is used to estimate the Hessian of the static map by multiplying it with the output $y(t)$ of the static map is defined in \cite{Ghaffari2011MultivariableNE} as
\begin{equation}
\hat{H}(t) = N(t)y(t)\;\; \text{with} \;\; N(t) = -\dfrac{8}{a^{2}}\cos{(2\omega t)}
\label{extrachapter.eq:hessian}
\end{equation}
whereas the signal $M(t)$ is used to estimate the gradient of the static map as follows:
\begin{equation}
    G(t) = M(t)y(t)\;\; \text{with} \;\; M(t) = \frac{2}{a}\sin{(\omega t)}.
\label{extrachapter.eq:ditherSM}
\end{equation}

\subsection{Additive Probing Signal}

The perturbation $S(t)$ is adapted from the basic ES to the case of PDE actuation dynamics. The trajectory generation problem as in \cite{Krsti2008BoundaryCO} is described as follows:
\begin{align}
    S(t) &\coloneqq -\partial_{x}\beta(0,t),\quad x \in (0,s(t)) \label{extrachapter.perturbation_beta}\\
    \partial_{t}\beta(x,t) &=          \partial_{xx}\beta(x,t)\label{extrachapter.stefan_beta}\\
    \beta(s(t),t) &= 0\label{extrachapter.boundary_beta}\\
    \beta_{x}(s(t),t) &= -a\omega\cos{(\omega t)}\label{extrachapter.initial_cond_beta},
\end{align}
where $\beta:[0,s(t)]\times\mathbb{R}_{+}\rightarrow\mathbb{R}$. The explicit solution of (\ref{extrachapter.perturbation_beta}) is found respectively for the reference trajectory and the reference solution postulated by a power series \cite{dunbar2003motion}:
\begin{align}
    &\xi(t) = a\sin{(\omega t)} \label{extrachapter.pertubation_s(t)}\\
    &\beta(x,t) = \sum_{i=0}^{\infty} \label{extrachapter.power_series} \dfrac{a_{i}(t)}{i!}[x-\xi(t)]^{i}.
\end{align}
We can calculate the first coefficients of the power series replacing the boundary conditions (\ref{extrachapter.boundary_beta}) and (\ref{extrachapter.initial_cond_beta}) at (\ref{extrachapter.power_series}), such that
\begin{equation}
    a_{0}(t) = 0,\;\;\;\;\; a_{1} = -\dot{\xi}(t).
\end{equation}
The general expression $a_{i}(t) = \dot{a}_{i-2}(t) - a_{i-1}(t)\dot{\xi}(t)$ is obtained by substituting (\ref{extrachapter.power_series}) at (\ref{extrachapter.stefan_beta}). We give the analytic expression of the first four coefficients of the series (\ref{extrachapter.power_series}) so that one can see how the successive derivatives of $\xi(t)$ appear
\begin{align}
    &a_{2}(t) = \dot{\xi}(t)^{2}\\
    &a_{3}(t) = \ddot{\xi}(t) - \dot{\xi}(t)^{3}\\
    &a_{4}(t) = \ddot{\xi}(t)^{2} + \ddot{\xi}(t)\dot{\xi}(t) + \dot{\xi}(t)^{4}.
\end{align}
The trajectory generation solution which provide all terms of the power series (\ref{extrachapter.power_series}) is given by \cite{hill1967parabolic}
\begin{equation}
    \beta(x,t) =  \sum_{i=0}^{\infty} \dfrac{1}{(2i)!}\dfrac{\partial^{i}}{\partial t^{i}}[x-\xi(t)]^{2i}. \label{extrachapter.beta_solution}
\end{equation}
Although (\ref{extrachapter.beta_solution}) is not an explicit expression, choosing suitable values for $a$ and $\omega$ in (\ref{extrachapter.pertubation_s(t)}), the series converges with few iterations of the infinite sum, getting the desirable sinusoid signal $\xi(t)$ in the output of the integrator.

According to (\ref{extrachapter.perturbation_beta}), we take the spatial derivative of (\ref{extrachapter.beta_solution}) and substitute $x=0$, therefore, we arrive at the final expression of
\begin{equation}
    S(t) = -\sum_{i=0}^{\infty} \dfrac{1}{(2i-1)!}\dfrac{\partial^{i}}{\partial t^{i}}[-a\sin{(\omega t)}]^{2i-1}.
\end{equation}

\subsection{Estimation Errors and PDE-Error Dynamics}

Since our objective is to find $\Theta^{*}$, which corresponds to the optimal unknown actuator $\theta(t)$, we introduce the following estimation errors

\begin{equation}
\hat{\theta}(t) = \theta(t) - S(t),\;\;\; \hat{\Theta}(t) = \Theta(t) - a\sin{(\omega t)},
\label{extrachapter.eq:estimated}
\end{equation}
\begin{equation}
\tilde{\theta}(t) \coloneqq \hat{\theta}(t) - \Theta^{*},\;\;\; \vartheta(t) \coloneqq \hat{\Theta}(t) - \Theta^{*},
\label{extrachapter.eq:estimated_errors}
\end{equation}
reminding that $\Theta(t) \coloneqq s(t)$. Combining $\hat{\Theta}(t)$ in (\ref{extrachapter.eq:estimated}) and (\ref{extrachapter.eq:estimated_errors}) we get the relation between the propagated estimation error $\vartheta(t)$, the propagated input $\Theta(t)$ and the optimizer of the static map $\Theta^{*}$

\begin{equation}
    \Theta(t) - \Theta^{*} = \vartheta(t) + a\sin{(\omega t)}. \label{extrachapter.vasin}
\end{equation}

Let us define

\begin{equation}
    u(x,t) = \alpha(x,t) - \beta(x,t) \label{extrachapter.u},
\end{equation}
\begin{equation}
    \hat{\theta}(t) = U(t). \label{extrachapter.utheta}
\end{equation}

By (\ref{extrachapter.Theta_p_actuator})-(\ref{extrachapter.theta_actuator}) and (\ref{extrachapter.perturbation_beta})-(\ref{extrachapter.initial_cond_beta}) with the help of (\ref{extrachapter.eq:estimated}) and (\ref{extrachapter.eq:estimated_errors}), we have our original system:
\begin{align}
    \dot{\vartheta}(t) &= -u_{x}(s(t),t),\quad x \in (0,s(t)) \label{extrachapter.dynamics_error}\\
    u_{t}(x,t) &= u_{xx}(x,t) \label{extrachapter.wavekv_error}\\
    u(s(t),t) &= 0\label{extrachapter.boundary_actuator}\\
    -u_{x}(0,t) &= U(t).\label{extrachapter.controller_error}
\end{align}

\section{Control Design} \label{extrachapter.controller}

\subsection{Stefan Compensation} \label{extrachapter.compensated}

We consider the PDE-ODE cascade (\ref{extrachapter.dynamics_error})-(\ref{extrachapter.controller_error}) and use the backstepping transformation 
\begin{equation}
\begin{split}
w(x,t) &= u(x,t) - \bar{K}\int_{x}^{s(t)} (x-\sigma)u(\sigma,t)\,dy\\
& -\bar{K}(x-s(t))\vartheta(t)
\label{extrachapter.eq:control_backstepping}
\end{split}
\end{equation}
with $\bar{K}>0$ is an arbitrary controller gain. Equation (\ref{extrachapter.eq:control_backstepping}) transforms (\ref{extrachapter.dynamics_error})-(\ref{extrachapter.controller_error}) into the target system:
\begin{align}
    \dot{\vartheta}(t) &= -\bar{K}\vartheta(t) - w_{x}(s(t),t),\,\,\,\,x \in (0,s(t)) \label{extrachapter.dynamics_target1}\\
    w_{t}(x,t) &= w_{xx}(x,t) + \bar{K}\dot{s}(t)\vartheta(t) \label{extrachapter.stefan_target}\\
    w_{x}(0,t) &= 0\label{extrachapter.boundary_target1}\\
    w(s(t),t) &= 0.\label{extrachapter.controller_target1}
\end{align}
The compensation controller can be obtained by taking the derivative of (\ref{extrachapter.eq:control_backstepping}) with respect to $t$ and $x$ respectively along of the solution (\ref{extrachapter.dynamics_error})-(\ref{extrachapter.controller_error}) and substituting $x=0$:
\begin{equation}
    U(t) = -\bar{K}\left( \vartheta(t) + \int_{0}^{s(t)} u(x,t)\,dx\right) \label{extrachapter.control_law}.
\end{equation}

\subsection{Implementable Extremum Seeking Control Law}
\label{coreancontrol}

Since we have no measurement on $\vartheta(t)$, (\ref{extrachapter.control_law}) is not applicable directly. Thus, introducing a result of \cite{Ghaffari2011MultivariableNE}, the average version of the gradient and Hessian estimates are calculated by

\begin{equation}
G_{\rm av}(t) = H\vartheta_{\rm av}(t), \;\;\; \hat{H}_{\rm av}(t) = H.
\label{extrachapter.eq:grad_hessian_estimate}
\end{equation}

Averaging (\ref{extrachapter.control_law}), choosing $\bar{K} = KH$ with $K<0$ and plugging the average gradient and Hessian estimates (\ref{extrachapter.eq:grad_hessian_estimate}), we obtain

\begin{equation}
    U_{\rm av}(t) = -KG_{\rm av}(t) - KH\int_{0}^{s_{\rm av}(t)} u_{\rm av}(x,t)\,dx. \label{extrachapter.average_control_law}
\end{equation}

We introduce a low-pass filter to the controller with the purpose of applying the average theorem for infinite-dimensional systems \cite{hale1990averaging} in the following stability proof, such that
\begin{equation}
U(t) = \dfrac{c}{s+c}\Bigg\{K\Bigg[G(t)+
\hat{H}(t) \int_{0}^{s(t)} u(x,t)\,dx \Bigg]\Bigg\},
\label{extrachapter.eq:filter_control_law}
\end{equation}
for $c>0$ sufficiently large.

\begin{figure}[!htb]
\begin{center}
    \includegraphics[width=8.5cm]{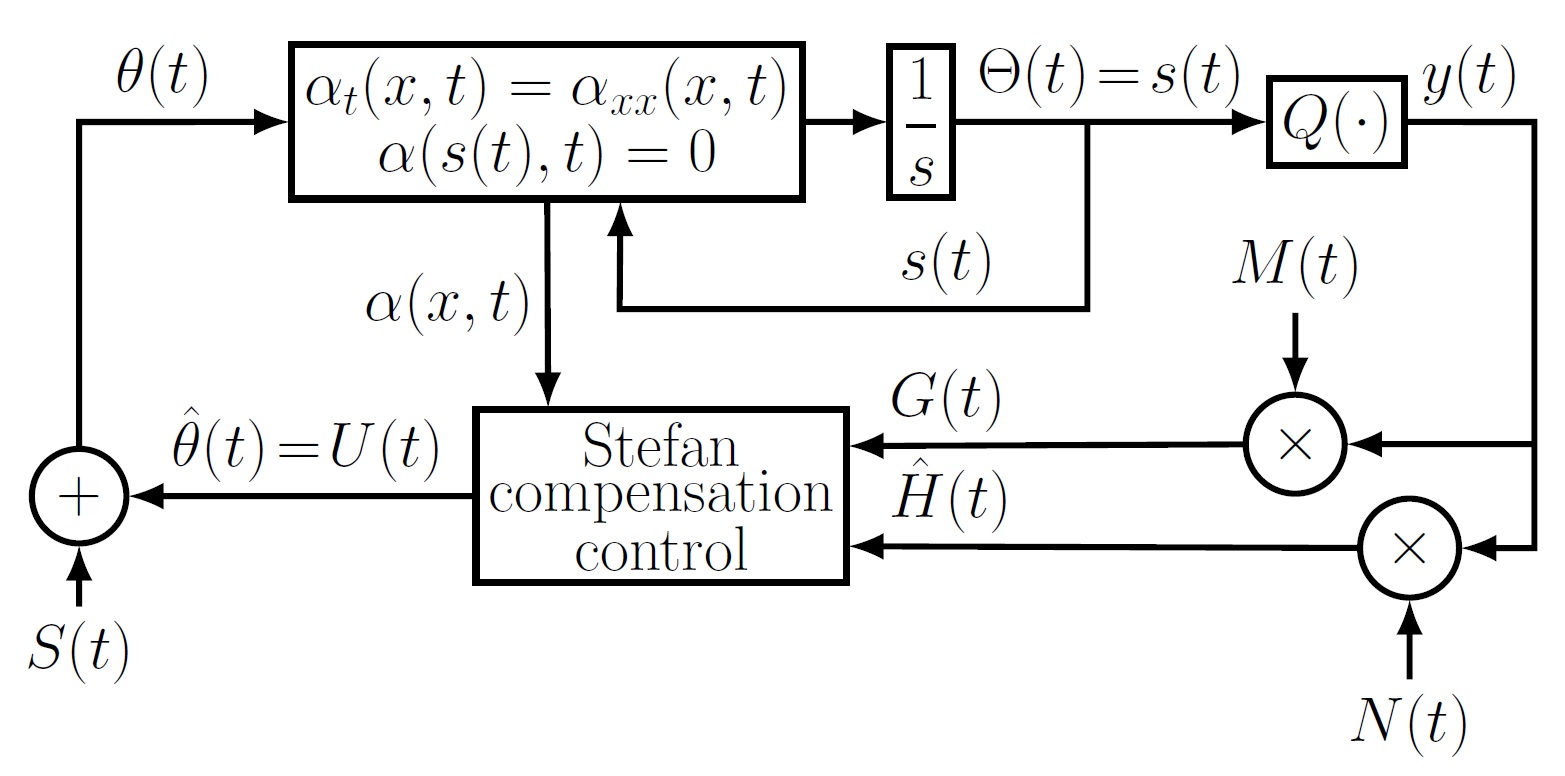} 
	\caption{ES control loop applied to the one-phase Stefan problem.}
	\label{extrachapter.fig:esc}
\end{center}
\end{figure}

Adapting the original scheme in \cite{Oliveira2017ExtremumSF} and combining (\ref{extrachapter.Theta_p_actuator})-(\ref{extrachapter.theta_actuator}) with the proposed boundary control law (\ref{extrachapter.eq:filter_control_law}), the closed-loop ES with actuation Stefan PDE dynamics is shown in Figure \ref{extrachapter.fig:esc}.

\section{Stability Analysis} \label{extrachapter.stability}

The following theorem summarizes the stability properties for the average version of the error-dynamics (\ref{extrachapter.dynamics_error})-(\ref{extrachapter.controller_error}).


\begin{thm} \label{reajeitando_teoremas}
Assume the model validity conditions $T_{\rm av}(x,t)>T_m$, $\dot{s}_{\rm av}(t)>0$, $x\in(0,s_{\rm av}(t))$ and
$s_0<s_{\rm av}(t)<s^*$ are satisfied at least in the average sense, $\forall t \geq 0$, $s^*=\Theta^*$ and for initial conditions $(T_{\rm av}(x,0),s_0)$ compatible with the control law $U(t)$ in (\ref{extrachapter.eq:filter_control_law}). Then, for a sufficiently large $c>0$, the average version of the closed-loop system
(\ref{extrachapter.dynamics_error})-(\ref{extrachapter.controller_error}) is exponentially stable in the sense of the norm $\| u _{\rm{av}}(t)\|^2_{\mathcal{H}_1} + \left | \vartheta_{\rm{av}}(t) \right |^{2}$, \textit{i.e.},  
\begin{equation}
    \|u_{\rm av}(t)\|_{\mathcal{H}_{1}}^{2} + \left|\vartheta_{\rm av}(t)\right|^{2} \leq M(\|u_{\rm av}(0)\|_{\mathcal{H}_{1}}^{2} + \vartheta_{\rm av}(0)^{2})e^{-nt}\,, \quad \forall t \geq 0\,, 
\end{equation}
and appropriate constants $M, n>0$. 
\end{thm}

\textbf{\textit{Proof.}}~ The proof is carried out in \textbf{Steps 1} to \textbf{3} below.

\subsubsection*{\textbf{Step 1:}\ \ \textit{Average closed-loop system}} 

The average version of the system (\ref{extrachapter.dynamics_error})-(\ref{extrachapter.controller_error}) is

\begin{align}
&\dot{\vartheta}_{\rm av}(t) = -(u_{\rm av})_{x}(s_{\rm av}(t),t),\,\,\,\, x \in (0,s_{\rm av}(t)) \label{extrachapter.dynamics_average_error}\\
&(u_{\rm av})_{t}(x,t) = (u_{\rm av})_{xx}(x,t) \label{extrachapter.stefan_average_error}\\
&u_{\rm av}(s_{\rm av}(t),t) = 0\label{extrachapter.boundary_average_error}\\
\begin{split}
&\dfrac{d}{dt}(u_{\rm av})_{x}(0,t) = -c(u_{\rm av})_{x}(0,t)\\
&- cKH\left[\vartheta_{\rm av}(t) + \int_{0}^{s_{\rm av}(t)} u_{\rm av}(x,t) \,dx \right],
\label{extrachapter.controller_average_error}
\end{split}
\end{align}
where the low-pass filter is represented in the state-space form. To derive (\ref{extrachapter.controller_average_error}), we plug the relationships $\vartheta(t)+a\sin(\omega t) = \Theta(t) - \Theta^{*}$, $G(t) = M(t)y(t)$ and (\ref{extrachapter.eq:hessian}) into (\ref{extrachapter.eq:filter_control_law}). With the help of the identities $2\sin^{2}(\omega t) = 1-\cos({2\omega t})$, $2\sin({\omega t})\cos({2\omega t}) = \sin({3\omega t}) - \sin({\omega t})$, $4\sin^{3}({\omega t}) = 3\sin({\omega t}) - \sin({3\omega t})$ and $4\sin^{2}(\omega t)\cos({2\omega t}) = 2\cos({2\omega t}) - \cos({4\omega t}) -1$ and applying averaging, we arrive at (\ref{extrachapter.controller_average_error}).\\

The backstepping transformation 

\begin{equation}
\begin{split}
w(x,t) &= u_{\rm av}(x,t) - KH\int_{x}^{s_{\rm av}(t)} (x-\sigma)u_{\rm av}(\sigma,t)\,d\sigma  \\
& -KH(x-s_{\rm av}(t))\vartheta_{\rm av}(t)
\label{extrachapter.eq:control_backstepping_avg}
\end{split}
\end{equation}
maps the average error-dynamics (\ref{extrachapter.dynamics_average_error})-(\ref{extrachapter.controller_average_error}) into the exponentially stable target system after assuming $c$ $\rightarrow$ +$\infty$ for the sake of simplicity. Consequently,

\begin{align}
    &\dot{\vartheta}_{\rm av}(t) = -KH\vartheta_{\rm av}(t) - w_{x}(s_{\rm av}(t),t),\,\,\,\, x \in (0,s_{\rm av}(t)) \label{extrachapter.dynamics_target2}\\
    &w_{t}(x,t) = w_{xx}(x,t) + KH\dot{s}_{\rm av}(t)\vartheta_{\rm av}(t) \label{extrachapter.stefan_target2}\\
    &w_{x}(0,t) = 0\label{extrachapter.boundary_target2}\\
    &w(s_{\rm av}(t),t) = 0.\label{extrachapter.controller_target2}
\end{align}

\subsubsection*{\textbf{Step 2:}\ \ \textit{Inverse transformation}} 

To ensure the equivalent stability property between the target system and the original system, the invertibility of the transformation (\ref{extrachapter.eq:control_backstepping_avg}) needs to be guaranteed. Assume the inverse transformation that maps (\ref{extrachapter.dynamics_target2})-(\ref{extrachapter.controller_target2}) into (\ref{extrachapter.dynamics_average_error})-(\ref{extrachapter.controller_average_error}):

\begin{equation}
\begin{split}
u_{\rm av}(x,t) &= w(x,t) + \int_{x}^{s_{\rm av}(t)} k(x-\sigma)w(\sigma,t)\,d\sigma  \\
& +\phi(x-s_{\rm av}(t))\vartheta_{\rm av}(t),
\label{extrachapter.eq:inverse_backstepping}
\end{split}
\end{equation}
where $k(x-\sigma)$ and $\phi(x-s_{\rm av}(t))$ are the kernel functions. Taking the derivatives with respect to $t$ and $x$, respectively, along the solution of (\ref{extrachapter.dynamics_target2})-(\ref{extrachapter.controller_target2}), the functions $\phi(x)$ and $k(x-\sigma)$ must satisfy:

\begin{equation}
    \phi^{''}(x) = -KH\phi(x),\,\,\, \phi(0) = 0,\,\,\, \phi^{'} = KH\,, \label{extrachapter.cond1}
\end{equation}
\begin{equation}
    k(x-s_{\rm av}(t)) = \phi(x-s_{\rm av}(t))\,, \label{extrachapter.cond2}
\end{equation}
\begin{equation}
    \phi^{'}(x-s_{\rm av}(t)) = KH\left(1+\int_{x}^{s_{\rm av}(t)} k(x-\sigma)\,d\sigma\right). \label{extrachapter.cond3}
\end{equation}

The solutions of the gain kernel can be deduced from (\ref{extrachapter.cond1})-(\ref{extrachapter.cond3}), such that

\begin{equation}
    \phi(x) = KH\sqrt{\dfrac{1}{KH}}\sin({\sqrt{KH}x})\,, \label{extrachapter.phi}
\end{equation}
\begin{equation}
    k(x-\sigma) = \phi(x-\sigma). \label{extrachapter.k}
\end{equation}

Hence, replacing (\ref{extrachapter.phi}) and (\ref{extrachapter.k}) into (\ref{extrachapter.eq:inverse_backstepping}), we have the following inverse transformation:
\begin{equation}
    \begin{split}
    &u_{\rm av}(x,t) = w(x,t) \\ &+\int_{x}^{s_{\rm av}(t)}KH\sqrt{\dfrac{1}{KH}}
    \sin({\sqrt{KH}(x-\sigma)})w(\sigma,t)\,d\sigma \\
    &+KH\sqrt{\dfrac{1}{KH}}\sin({\sqrt{KH}(x-s_{\rm av}(t))})\vartheta_{\rm av}(t).
    \end{split}
\end{equation}

\subsubsection*{\textbf{Step 3:}\ \ \textit{Exponential stability}} 

We prove the exponential stability of the average closed-loop system based on the target system (\ref{extrachapter.dynamics_average_error})-(\ref{extrachapter.controller_average_error}) using the Lyapunov method. We consider the following Lyapunov functional:

\begin{equation}
    V = V_{1} + V_{2} + V_{3}\,, \label{extrachapter.V}
\end{equation}
\begin{equation}
    V_{1} = \dfrac{1}{2}\int_{0}^{s_{\rm av}(t)}w(x,t)^{2}\,dx\,, \label{extrachapter.V1}
\end{equation}
\begin{equation}
    V_{2} = \dfrac{1}{2}\int_{0}^{s_{\rm av}(t)}w_{x}(x,t)^{2}\,dx\,, \label{extrachapter.V2}
\end{equation}
\begin{equation}
    V_{3} = \rho\dfrac{1}{2}\vartheta_{\rm av}(t)^{2}\,. \label{extrachapter.V3}
\end{equation}

Taking the derivative of (\ref{extrachapter.V1}) with respect to $t$:

\begin{equation}
    \begin{split}
    \dot{V_{1}} &= -\int_{0}^{s_{\rm av}(t)}w_{x}(x,t)^{2}\,dx \\ &+KH\dot{s}_{\rm av}(t)\vartheta_{\rm av}(t)\int_{0}^{s_{\rm av}(t)}w(x,t)\,dx. \label{extrachapter.dot_V1}
    \end{split}
\end{equation}

Taking the derivative of (\ref{extrachapter.V2}) with respect to $t$:

\begin{equation}
    \begin{split}
    &\dot{V_{2}} = w_{x}(s_{\rm av}(t),t)w_{t}(x,t)+\dfrac{1}{2}\dot{s}_{\rm av}(t)w_{x}(s_{\rm av}(t),t)^{2}  \\ &\!-\!KH\dot{s}_{\rm av}(t)\vartheta_{\rm av}(t)w_{x}(s_{\rm av}(t),t) - \int\displaylimits_{0}^{s_{\rm av}(t)}w_{xx}(x,t)^{2}\,dx . \label{extrachapter.half_dot_V2}
    \end{split}
\end{equation}

Using the relationship $w_{t}(s_{\rm av}(t),t) = -\dot{s}_{\rm av}(t)w_{x}(s_{\rm av}(t),t)$ and replacing it into (\ref{extrachapter.half_dot_V2}), we obtain
\begin{equation}
    \begin{split}
    \dot{V_{2}} &= -\int_{0}^{s_{\rm av}(t)}w_{xx}(x,t)^{2}\,dx -\dfrac{1}{2}\dot{s}_{\rm av}(t)w_{x}(s_{\rm av}(t),t)^{2} \\
    &-KH\dot{s}_{\rm av}(t)\vartheta_{\rm av}(t)w_{x}(s_{\rm av}(t),t). \label{extrachapter.dot_V2}
    \end{split}
\end{equation}

Taking the derivative of (\ref{extrachapter.V3}) with respect to $t$, lead us to

\begin{equation}
    \dot{V_{3}} = -\rho KH\vartheta_{\rm av}(t)^{2} - \rho \vartheta_{\rm av}(t)w_{x}(s_{\rm av}(t),t). \label{extrachapter.dot_V3}
\end{equation}

Substituting the terms (\ref{extrachapter.dot_V1}), (\ref{extrachapter.dot_V2}) and (\ref{extrachapter.dot_V3}) into the time derivative of (\ref{extrachapter.V}) and using the Young's inequality in $-\rho \vartheta_{\rm av}(t)w_{x}(s_{\rm av}(t),t)$, $KH\dot{s}_{\rm av}(t)\vartheta_{\rm av}(t)\int_{0}^{s_{\rm av}(t)}w(x,t)\,dx$ and $-KH\dot{s}_{\rm av}(t)\vartheta_{\rm av}(t)w_{x}(s_{\rm av}(t),t)$, we have

\begin{equation}
    \begin{split}
    \dot{V} &\leq -\int_{0}^{s_{\rm av}(t)}w_{xx}(x,t)^{2}\,dx - \int_{0}^{s_{\rm av}(t)}w_{x}(x,t)^{2}\,dx\\
    &- \dfrac{\rho KH}{2}\vartheta_{\rm av}(t)^{2} + \dfrac{\rho}{2KH}w_{x}(s_{\rm av}(t),t)^{2}\\
    &+ \dot{s}_{\rm av}(t)\left(\dfrac{s^{*}}{2}\int_{0}^{s_{\rm av}(t)}w(x,t)^{2} + (KH)^{2}\vartheta_{\rm av}(t)^{2}\right). \label{extrachapter.half_dot_V}
    \end{split}
\end{equation}
By choosing $\rho = \dfrac{KH}{4s^{*}}$ and applying Poincar\'{e} and Agmon's inequalities at $\int_{0}^{s_{\rm av}(t)}w(x,t)^{2}\,dx$ and $w_{x}(s_{\rm av}(t),t)^{2}$, respectively, we obtain
\begin{equation}
    \begin{split}
    \dot{V} \leq& -\dfrac{1}{8{s^{*}}^{2}}\int\displaylimits_{0}^{s_{\rm av}(t)}w_{x}(x,t)^{2}dx- \dfrac{1}{4{s^{*}}^{2}}\int\displaylimits_{0}^{s_{\rm av}(t)}w(x,t)^{2}dx  \\
    &+\dot{s}_{\rm av}(t)\left(\dfrac{s^{*}}{2}\int_{0}^{s_{\rm av}(t)}\!\!\!w(x,t)^{2}\,dx \!+\! (KH)^{2}\vartheta_{\rm av}(t)^{2}\right) - \dfrac{\rho KH}{2}\vartheta_{\rm av}(t)^{2} \\
    \leq& -mV + n\dot{s}_{\rm av}(t)V, \label{extrachapter.dot_V}
    \end{split}
\end{equation}
where
\begin{equation}
    n = \max\,\, \{1,\,8s^{*}KH\},\,\,\,\,\,m = \min\, \{1/4{s^{*}}^{2},\,KH\}.
\end{equation}
The term $n\dot{s}_{\rm av}(t)V$ on the right-hand side of (\ref{extrachapter.dot_V}) does not let us to directly conclude exponential stability. To deal with it, a new Lyapunov function candidate $W$ is defined according to
\begin{equation}
    W(t) = V(t)e^{-ns_{\rm av}(t)}. \label{extrachapter.W}
\end{equation}
The time derivative of (\ref{extrachapter.W}) can be calculated using (\ref{extrachapter.dot_V}) such that 
\begin{equation}
    \dot{W}(t) = (\dot{V}(t) - n\dot{s}_{\rm av}(t)V(t))e^{-ns_{\rm av}(t)} \leq -m W(t).
\end{equation}
Taking into account (\ref{extrachapter.V}), we can establish the following relationship:
\begin{equation}
    \|w(t)\|_{\mathcal{H}_{1}}^{2} + \rho \vartheta_{\rm av}(t)^{2} \!\leq\! e^{ns^{*}}(\|w_{0}\|_{\mathcal{H}_{1}}^{2} + \rho \vartheta_{\rm av}(0)^{2})e^{-mt}\,, \quad w(x,0)\!=\!w_{0}\,. \label{extrachapter.norm_w}
\end{equation}
Hence, we can conclude the existence of a positive constant $M>0$ using the inverse transformation (\ref{extrachapter.eq:control_backstepping_avg}) combined with Young's and Cauchy-Schwarz inequalities, such that
\begin{equation}
    \|u_{\rm av}(t)\|_{\mathcal{H}_{1}}^{2} + \vartheta_{\rm av}(t)^{2} \!\leq\! M(\|u_{\rm av0}\|_{\mathcal{H}_{1}}^{2} \!+\! \vartheta_{\rm av}(0)^{2})e^{-mt}\,, \quad \!u_{\rm{av}}(x,0)\!=\!u_{\rm av0}\,, \label{extrachapter.norm_u}
\end{equation}
which completes the proof. \hfill $\square$

\section{Asymptotic Convergence to a Neighborhood of the\\ Extremum Point} \label{samsungs23+}

In Theorem~\ref{reajeitando_teoremas}, we prove that the average closed-loop
system (\ref{extrachapter.dynamics_average_error})-(\ref{extrachapter.controller_average_error}) is exponentially stable. However, there is no suitable averaging theorem for moving-boundary PDE systems and it remains as an open problem in the literature. 


If this theorem existed such that employed for PDEs of fixed domains \cite{hale1990averaging}, then we would apply it to (\ref{extrachapter.dynamics_error})-(\ref{extrachapter.controller_error}) and conclude for the non-average system the existence of a unique exponentially stable periodic solution in $t$ of period $\Pi \coloneqq 2\pi/\omega$, denoted by $\vartheta^{\Pi}(t)$, $u^{\Pi}(x,t)$, satisfying 
\begin{equation}
\begin{aligned}
& {} \left(|\vartheta^{\Pi}(t)|^{2}+\|u^{\Pi}(t)\|^{2}+\|u_{x}^{\Pi}(t)\|^{2} \right)^{1/2} \leq \mathcal{O}(1/\omega)\,, \quad \forall t \geq 0\,.
\label{extrachapter.eq:order_period}
\end{aligned}
\end{equation}
%
%
%
On the other hand, the asymptotic convergence to a neighborhood of the extremum point would be proved taking the absolute value of (\ref{extrachapter.vasin}):
\begin{equation}
    |\Theta(t)-\Theta^{*}| = |\vartheta(t) + a\sin{(\omega t)}|, \label{extrachapter.vasin2}
\end{equation} 
and writing (\ref{extrachapter.vasin2}) in terms of the periodic solution $\vartheta^{\Pi}(t)$: $|\Theta(t)\!-\!\Theta^{*}| \!=\! |\vartheta(t)\!-\!\vartheta^{\Pi}(t)\!+\!\vartheta^{\Pi}(t) \!+\! a\sin{(\omega t)}|$. 
By applying 
again the appropriate averaging theorem, one would have $\vartheta(t)\!-\!\vartheta^\Pi(t)\!\to\!0$ exponentially and, consequently, 
\begin{equation}
    \limsup_{t \to \infty} |\Theta(t)\!-\!\Theta^{*}| = \limsup_{t \to \infty} |\vartheta^{\Pi}(t) + a\sin{(\omega t)}| \label{extrachapter.Order_Theta(t)1}.
\end{equation}
Finally, with (\ref{extrachapter.eq:order_period}) we would arrive to 
\begin{align}
\limsup_{t \to \infty} |\Theta(t)-\Theta^{*}| = \mathcal{O}\left(|a|+1/\omega\right). \label{extrachapter.Order_Theta(t)}
\end{align}
In order to show the convergence of the output $y(t)$, we could follow the same steps employed for $\Theta(t)$ by plugging (\ref{extrachapter.vasin}) into   (\ref{extrachapter.eq:final_output_static_map}), such that
\begin{equation}
\begin{split}
\limsup_{t \to \infty} |y(t)-y^{*}| = \limsup_{t \to \infty} |H\vartheta^{2}(t) + Ha^{2}\sin{(\omega t)}^{2}|. \label{extrachapter.Order_y(t)1}
\end{split}
\end{equation}
Hence, rewriting (\ref{extrachapter.Order_y(t)1}) in terms of $\vartheta^{\Pi}(t)$ and again with the help of (\ref{extrachapter.eq:order_period}), we finally get 
\begin{align}
\limsup_{t \to \infty} |y(t)-y^{*}| = \mathcal{O}\left(|a|^{2}+1/\omega^{2}\right). \label{extrachapter.order_y(t)}
\end{align}

\section{Delay-Compensated Control} \label{delay_coreano}

In this section, we formulate the Stefan problem for ES \textit{with delays} in the actuation dynamics. We have omitted it before since the formulation of the Stefan problem
for ES \textit{without} delays itself has not been done before in the literature yet either. Moreover, the delay would complicate the understanding of the new contribution for the Stefan model. In this delayed scenario, we show in the following the relation between the designed control law and a state prediction used for delay compensation.

\begin{figure}[htb!]
\begin{center}
\includegraphics[width=8cm]{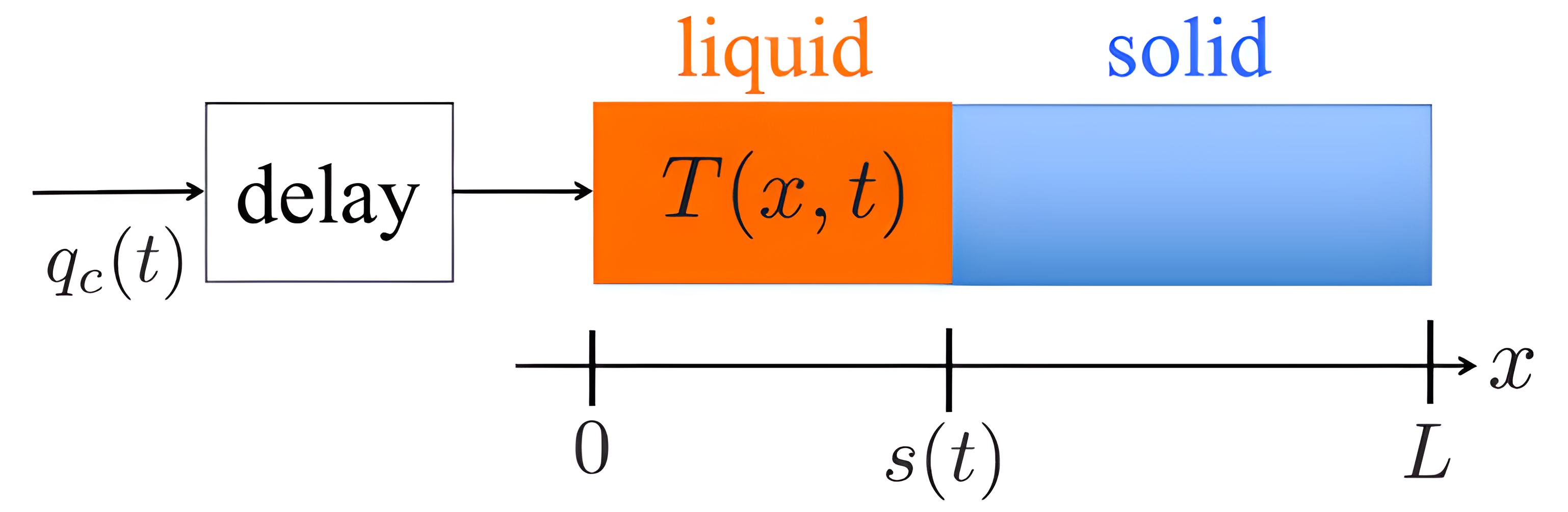}    
\caption{Schematic of one-phase Stefan problem with actuator delay \citep{koga2020delay}.}
\label{fig:stefan_delay}
\end{center}
\end{figure}

The diffusion equation of the temperature in the liquid-phase with actuator delay represented by Figure \ref{fig:stefan_delay} is described by:

\begin{align}
    T_{t}(x,t) &= \alpha T_{xx}(x,t),\quad x\in (0,s(t)),\quad \alpha = \dfrac{k}{\rho C_{p}}\label{T1_delay}\\
    -kT_{x}(0,t) &= q_{c}(t-D)\label{T2_delay}\\
    T(s(t),t) &= T_{m}\label{T3_delay}\\
    \dot{s}(t) &= -\beta T_{x}(s(t),t),\quad \beta = \dfrac{k}{\rho \Delta H^{*}}, \label{T4_delay}
\end{align}
where $D$ is the input delay.

As assumed in Section \ref{coreancontrol},  $\alpha = \beta = k = 1$. Based on reference \citep{koga2020delay}, the control law $U(t)$ for our extremum seeking problem of the closed-loop system consisting of the plant (\ref{T1_delay})-(\ref{T4_delay}) would be:

\begin{equation}
\begin{split}
     U(t) = -\bar{K} \bigg(\int_{t-D}^{t}U(\psi)\,d\psi +\int_{0}^{s(t)}u(x,t)\,dx + \vartheta(t)\bigg). \label{U_delay}
\end{split}
\end{equation}

In the next, we are going to represent equation (\ref{extrachapter.control_law}) as $U_{d}(t)$ to differ from (\ref{U_delay}). The objective of this section is to prove the equivalence between the delay-compensated control (\ref{U_delay}) and the predictor-based feedback. In other words, $U(t) \equiv U_{d}(t+D)$, which can be rewritten as

\begin{equation}
     U_{d}(t+D) = -\bar{K} \left(\int\displaylimits_{0}^{s(t+D)}u(x,t+D)\,dx+\vartheta(t+D)\right). \label{U_predictor}
\end{equation}

Integrating $\dot{s}(t) = -u_{x}(s(t),t)$, from $t$ to $t+D$, yields

\begin{equation}
     s(t+D) = s(t) -\int_{t}^{t+D}u_{x}(s(\tau),\tau)\,d\tau. \label{s_predictor}
\end{equation}

Now, integrating $u_{t}(x,t) = u_{tt}(x,t)$ in time, from $t$ to $t+D$, and in space, from $0$ to $s(t+D)$, respectively, we get

\begin{equation}
\begin{split}
     \int_{0}^{s(t+D)}&\int_{t}^{t+D}u_{t}(x,t)\,dt\,dx = \\ &\int_{0}^{s(t+D)}\int_{t}^{t+D}u_{xx}(x,t)\,dt\,dx \label{U_predictor1}
\end{split}
\end{equation}
\begin{equation}
\begin{split}
     &\int_{0}^{s(t+D)}u(x,t+D)\,dx = \int_{0}^{s(t+D)}u(x,t)\,dx\\ &+\int_{t}^{t+D}u_{x}(s(t+D),\tau)\,d\tau + \int_{t-D}^{t}U_{d}(\varphi)\,d\varphi. \label{U_predictor2}
\end{split}
\end{equation}

Substituting (\ref{s_predictor}) and (\ref{U_predictor2}) in (\ref{U_predictor}), we obtain

\begin{equation}
\begin{split}
     &U_{d}(t+D) = -\bar{K}\bigg(\int_{0}^{s(t+D)}u(x,t)\,dx \\ &+\int_{t}^{t+D}(u_{x}(s(t+D),\tau) - u_{x}(s(\tau),\tau))\,d\tau\\ &+ \int_{t-D}^{t}U_{d}(\varphi)\,d\varphi + \vartheta(t)\bigg). \label{U_predictor3}
\end{split}
\end{equation}

Analyzing the second integral of the right-hand side in (\ref{U_predictor3}):
\begin{equation}
\begin{split}
     &\int_{t}^{t+D}(u_{x}(s(t+D),\tau) - u_{x}(s(\tau),\tau))\,d\tau = \\ &\int_{t}^{t+D}\int_{s(\tau)}^{s(t+D)}u_{xx}(x,\tau)\,dx\,d\tau = \\
     & \int_{s(t)}^{s(t+D)}(u(x,s^{-1}(x))-u(x,t))\,dx.\label{U_predictor4}
\end{split}
\end{equation}

The boundary condition $u(s(t),t) = 0,\, \forall t \geq 0$, implies that $u(x,s^{-1}(x)) = 0$. By applying this condition to (\ref{U_predictor4}), one has

\begin{equation}
\begin{split}
     &\int_{t}^{t+D}(u_{x}(s(t+D),\tau) - u_{x}(s(\tau),\tau))\,d\tau = \\ &-\int_{s(t)}^{s(t+D)}u(x,t)\,dx.\label{U_predictor5}
\end{split}
\end{equation}

Finally, replacing (\ref{U_predictor5}) into (\ref{U_predictor3}), we arrive at

\begin{equation}
     U_{d}(t) = -\bar{K} \left(\int\displaylimits_{t-D}^{t}U(\psi)\,d\psi +\int\displaylimits_{0}^{s(t)}u(x,t)\,dx + \vartheta(t)\right), \label{U_predictor6}
\end{equation}
which is the same as the delay-compensated control (\ref{U_delay}).

\begin{figure}[htb!]
\begin{center}
\includegraphics[width=9cm]{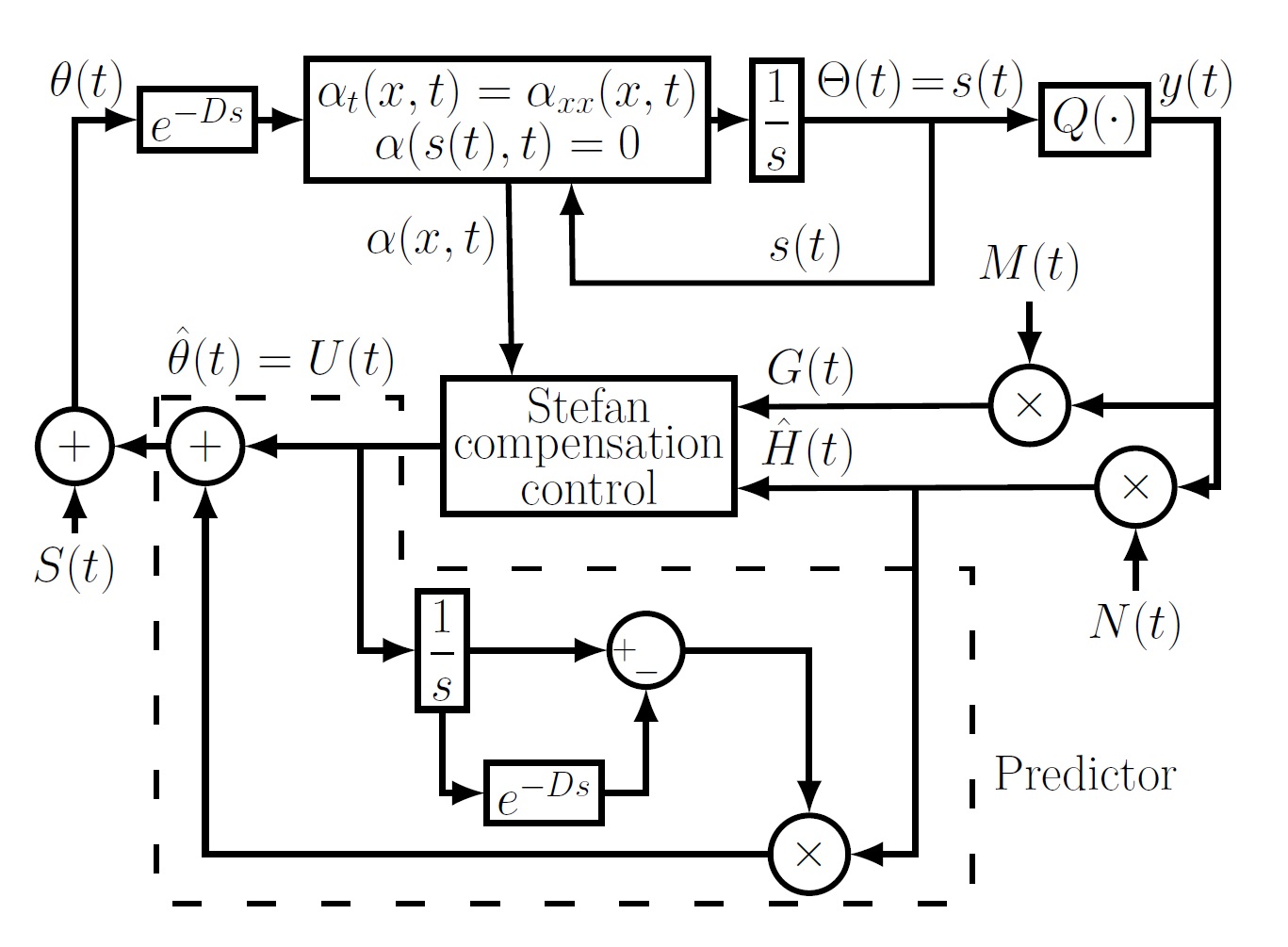}
\caption{Extremum seeking control loop with basic prediction scheme applied to the one-phase Stefan problem under actuation delay.} 
\label{fig:esc_delay}
\end{center}
\end{figure}

Comparing the prediction control law (\ref{U_predictor6}) with our original control law without delay (\ref{extrachapter.control_law}), we can verify that the additional expression on the right-hand side $\int_{t-D}^{t}U(\psi)\,d\psi$ is similar to the integral of \citep{Oliveira2017ExtremumSF} [Eq. 35] and it can be represented as a basic prediction scheme according to Figure \ref{fig:esc_delay}.


The dither signal $S(t)$ in Figure \ref{fig:esc_delay}, using as reference (\ref{extrachapter.beta_solution}), will be simply an implementable and advanced version of (\ref{extrachapter.perturbation_beta}), \textit{i.e.},

\begin{equation}
    S(t) = \beta_{x}(0,t+D).
\end{equation}

Using the same steps of Section \ref{extrachapter.controller} for (\ref{U_delay}), yields to the following implementable control law with a low-pass filter:

\begin{equation}
\begin{split}
U(t) &= \dfrac{c}{s+c}\Bigg\{K\Bigg[G(t)+
\hat{H}(t) \int_{t-D}^{t} U(\psi)\,d\psi \\
&+\hat{H}(t) \int_{0}^{s(t)} u(x,t)\,dx\Bigg]\Bigg\}.
\label{eq:filter_control_law_delay}
\end{split}
\end{equation}

\section{Simulations} \label{extrachapter.simulation}

The numerical simulation employs the quadratic map described in (\ref{extrachapter.eq:static_map}) and the parameters are chosen as stated in Table \ref{extrachapter.tb:margins}. For the sake of completeness, we have considered the general case containing a delay $D=0.5$ seconds in the actuation dynamics. It is worth to mention that this amount of delay is long enough to destabilized the closed-loop system, if not properly compensated.

\begin{table}[hb]
\begin{center}
\caption{Simulation parameters}\label{extrachapter.tb:margins}
\begin{tabular}{cccc}
 & \!\!\!Symbol & \!\!\!\!\!\!\!Description & Value \\\hline
                      & \!\!\!$K$ & \!\!\!\!\!\!\!controller gain & -0.1 \\
       Controller               & \!\!\!$a$ & \!\!\!\!\!\!\!perturbation amplitude & 0.1\\
     parameters       & \!\!\!$c$ & \!\!\!\!\!\!\!pole of the low-pass filter [rad/s] & 10\\                  & \!\!\!$\omega$ & \!\!\!\!\!\!\!perturbation frequency                     [rad/s] & 10\\
     
                      \hline
                      & \!\!\!$L$ & \!\!\!\!\!\!\!spatial domain & 1 \\ 
                & \!\!\!$\Theta^{*}$ & \!\!\!\!\!\!\!optimizer static map & 0.8 \\
    System         & \!\!\!$y^{*}$ & \!\!\!\!\!\!\!optimal value static map & 4 \\
    parameters                  & \!\!\!$H$ & \!\!\!\!\!\!\!Hessian of the static map & -1 \\
                      & \!\!\!$s_{0}$ & \!\!\!\!\!\!\!Initial interface [m] & 0.12 \\
                      & \!\!\!$T_{0}$ & \!\!\!\!\!\!\!Initial temperature [$^{\circ}C$] & 110 \\
                      & \!\!\!$T_{m}$ & \!\!\!\!\!\!\!Melting temperature [$^{\circ}C$] & 100\\
                      & \!\!\!$D$ & \!\!\!\!\!\!\!delay [s] & 0.5
\end{tabular}
\end{center}
\end{table}

\begin{figure*}[htb!]
\begin{center}
\includegraphics[width=14cm]{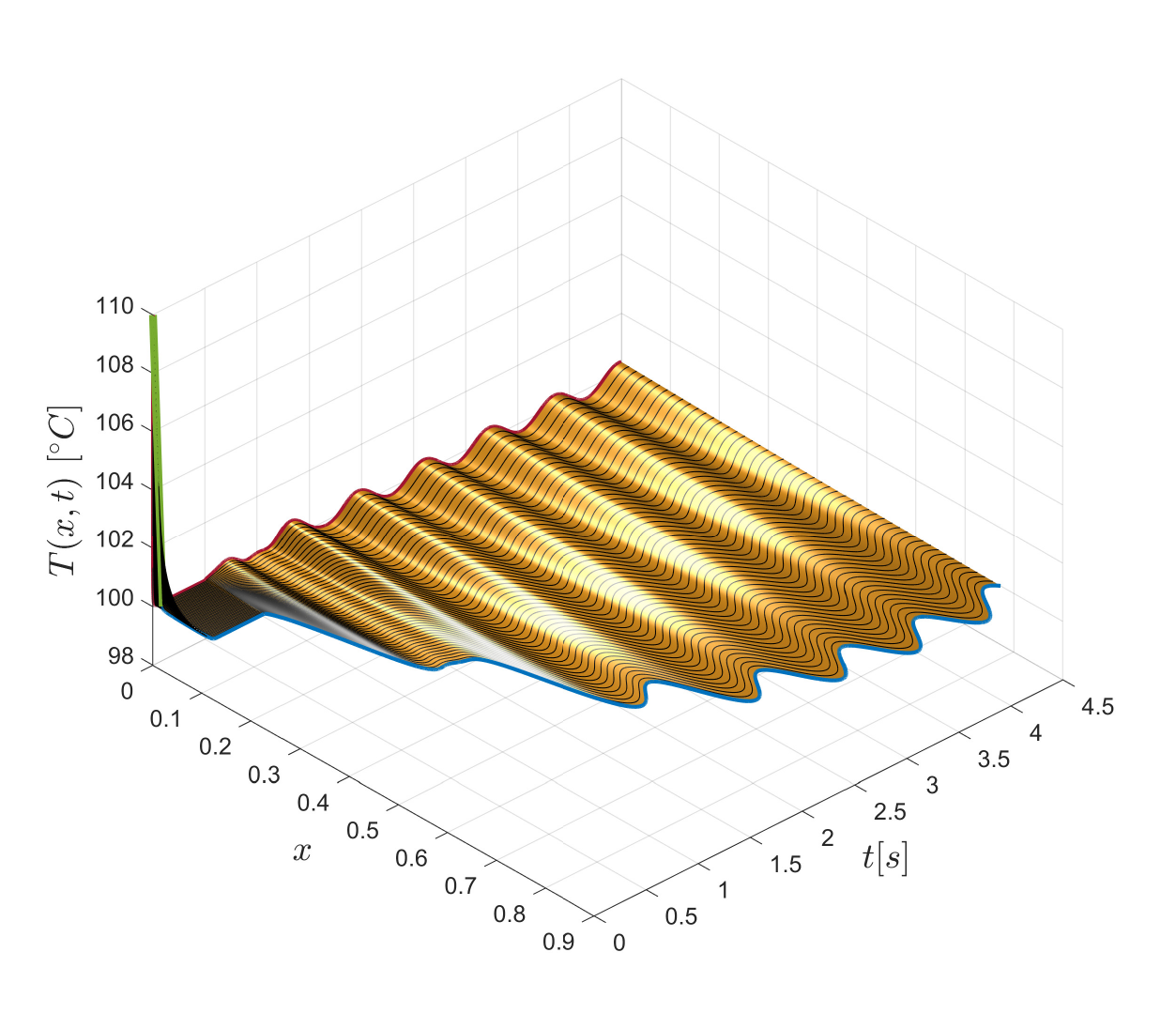} 
\vspace{-1cm}
\caption{The initial states $T(x,0)$ (green curve), $T(0,t)$ (red curve) and the convergence of $T(s(t),t)$ (blue curve) to $T_{m}=100^{\circ}C$ in a
three-dimensional space for the PDE state $T(x,t)$. The blue curve shows the expansion of the domain of the Stefan PDE. }
\label{extrachapter.fig:Theta_freedelay}
\end{center}
\end{figure*}

Figure \ref{extrachapter.fig:Theta_freedelay} corresponds to the numerical plot of the temperature profile for the closed-loop system converging in a three-dimensional space (taking into account the domain $L$ and the time $t$) to a close neighborhood of $T_{m}$.

\begin{figure*}[htb!]
\begin{center}
\includegraphics[width=14cm]{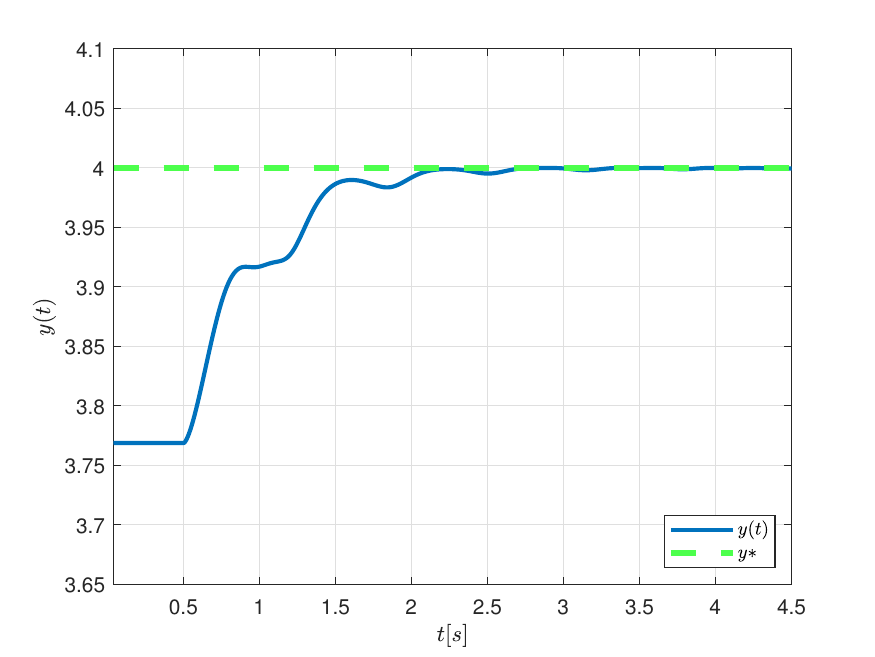}
\caption{Convergence of the output $y(t)$ to $y^{*}$.}
\label{extrachapter.fig:y_t}
\end{center}
\end{figure*}

\begin{figure*}[htb!]
\begin{center}
\includegraphics[width=14cm]{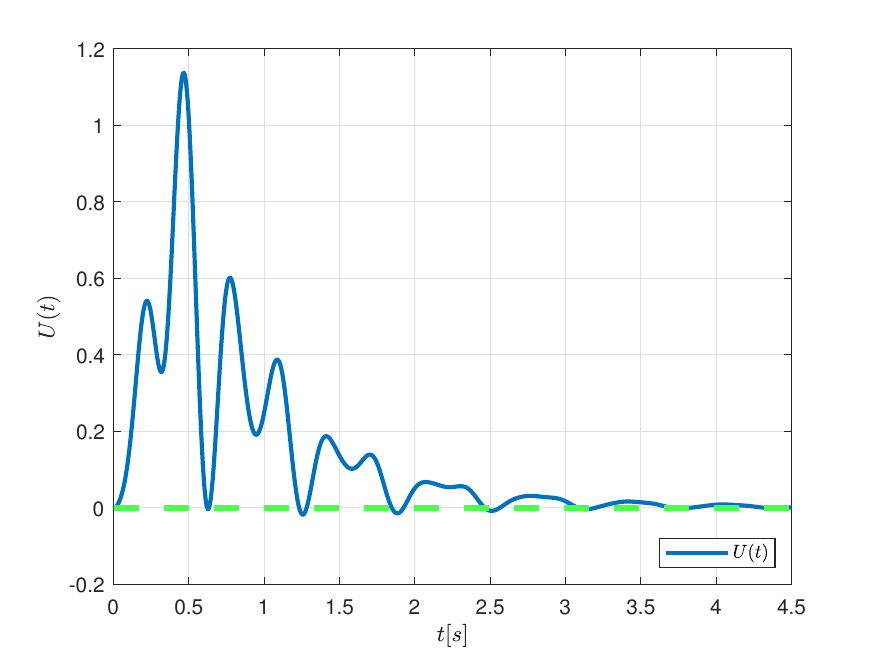} 
\caption{Convergence of the control signal $U(t)$ to $0$.} 
\label{extrachapter.fig:u_t}
\end{center}
\end{figure*}

Figure \ref{extrachapter.fig:u(s(t),t)} shows the convergence of the moving boundary to the optimizer $\Theta^{*}$. 
The sinusoidal movement of $s(t)$ would violate the usual conditions for the Stefan problem that the temperature remains above or below the melting temperature on the whole interval $[0, s(t)]$, forming a periodic chain of liquid and solid. However, it is well known in physics that the phase transition is not sharp but gradual. For a small perturbation amplitude, the system will be operating in the ``phase transition range'' when it reaches the periodic steady state at the extremum. A more detailed model is needed to capture the system dynamics in that narrow phase transition range. 
%
%
In addition, we could redesign the algorithm in order to introduce vanishing probing signals and tapering off the perturbation after the extremum neighborhood is achieved, as studied in \cite{durr2013lie}, \cite{scheinker2014non} and \cite{wang2016stability}. 
At last, Figures \ref{extrachapter.fig:y_t} and \ref{extrachapter.fig:u_t} show the convergence of the output $y(t)$ to $y^{*}$ and $U(t)$ to 0, respectively.

\begin{figure*}[htb!]
\begin{center}
\includegraphics[width=12cm]{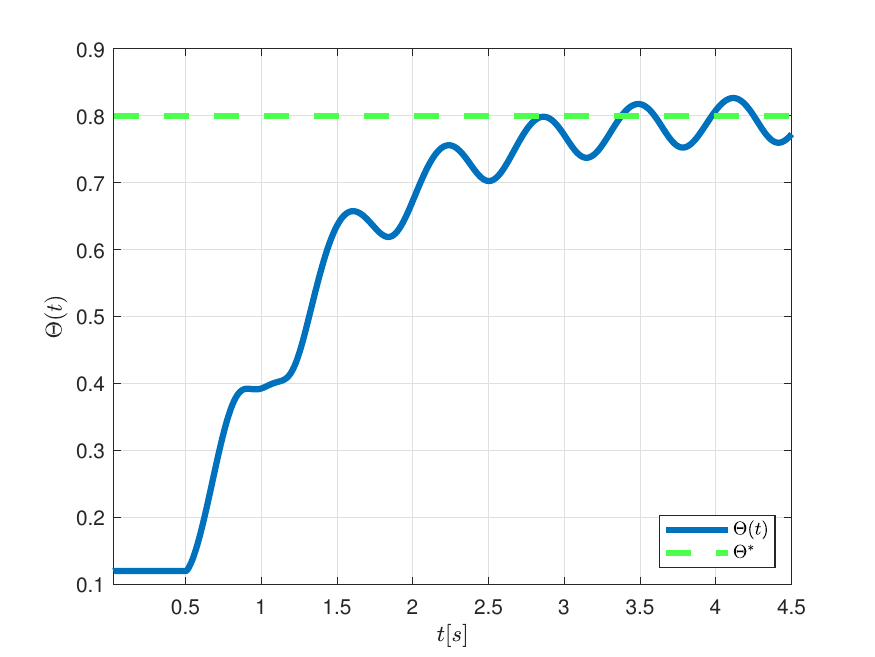} 
\caption{Convergence of $s(t):=\Theta(t)$ to a small neighborhood of $\Theta^*$.} 
\label{extrachapter.fig:u(s(t),t)}
\end{center}
\end{figure*}



\section{Conclusions} \label{concludingfinally}

The proposed approach manages to maximize the static map by searching the extremum point even in the presence of a Stefan PDE with moving boundary, possibly
including actuator delays as well. Although the actuation dynamics must be known, no information is assumed from the map parameters. The average boundary control law to compensate the actuation dynamics employed the backstepping methodology. Local exponential stability of the average system was guaranteed and convergence to a small neighborhood of the extremum was also indicated. An illustration of the benefits of the new extremum seeking scheme for the Stefan PDE is presented using consistent simulation results.

                                                   






  


\end{document}